\documentclass{amsart}
\usepackage{amssymb, array, verbatim, amscd}
\usepackage{amsmath, amsfonts,amsthm}

\newtheorem{theorem}{Theorem}[section]

\newtheorem{definition}[theorem]{Definition}
\newtheorem{example}[theorem]{Example}

\newcommand{\Z}{{\mathbb Z}}
\newcommand{\R}{{\mathbb R}}
\begin{document}

\title{Special homomorphisms between Probabilistic Gene Regulatory Networks }

\author{Mar\'{\i}a Alicia Avi\~n\'o }
\address{ Department of
Mathematic-Physics,
       University of Puerto Rico,
       Cayey, PR 00736}
\thanks{This research was supported by the National Institute of Health,
 PROGRAM SCORE, 2004-08, 546112, University of Puerto Rico-Rio Piedras Campus, IDEA Network of Biomedical
  Research Excellence, and the Laboratory Gauss University of Puerto Rico Research.
   I want to thank Professor E. Dougherty for his useful suggestions, and Professor O. Moreno for  his support
   during the last four years. }
%\ninept
\email{
mavino@cayey.upr.edu}
 \date{\today}

 \begin{abstract} In this paper we study finite dynamical systems with $n$
 functions acting on the same set $X$,  and  probabilities assigned to these functions, that it is called Probabilistic Regulatory Gene Networks (PRN) in \cite{D2005}.
 This concept is  the same or a natural generalization of the  concept
  Probabilistic Boolean Networks (PBN), introduced by I. Shmulevich, E. Dougherty, and W. Zhang in \cite{SDZ},   particularly  the model PBN has been using to describe genetic
 networks and has therapeutic applications, see \cite{SGHDZ}.
In PRNs  the most important question is to describe the steady
states of the systems, so in this paper we pay attention to the
idea of transforming a network to another without lost all the
properties,  in particular the probability distribution. Following
this objective we develop the concepts of homomorphism and
$\epsilon$-homomorphism of probabilistic regulatory networks,
 since these concepts bring the properties from
one networks to another.Projections are special homomorphisms, and hey always induce
invariant subnetworks that contain all cycles and steady  states in the network
.
\end{abstract}
\maketitle \section*{Introduction}
Genes can be understanding in their complexity behavior using models according with their discrete or continuous action. Developing computational tools permits describe gene functions and understand the mechanism of regulation \cite{48,49}. \emph{This understanding will have a significant impact on the development of techniques for drugs testing and therapeutic intervention  for treating human diseases}\cite{D2005,50,SGHDZ}.

We focus our attention in the discrete structure of  genetic regulatory networks, instead of, its dual moving continuo-discrete.
Probabilistic Gene Regulatory Network(PRgN)  is  a natural generalizations of the model Probabilistic Boolean Network (PBN), introduced by I. Shmulevich, E. Dougherty, and W. Zhang in \cite{SDZ}. The mathematical background of the model PgRN, is introduced here, for simplicity we work with functions defined over  a set $X$ to itself, with probabilities assigned to these functions. $X$ is a set of states of genes, for example $X=\{0,1\}^n$, if our network is a Boolean network. Working in this way, we can observe the dynamic of the network indeed focus our attention in the description of functions. The set $X$ can be a subset of $\{0,1\}^n$, and we can extend some classical ideas to regulatory network, such as invariant subnetworks, automorphisms group, etc.  In particular if $X$  is a vector space over a finite field, the functions are lineal functions, then we can use linear Algebra to describe the state space. Mapping are important in the study of networks, because they permit to recognize subnetworks, in particular determine when two networks are similar or equivalent. Special mappings are homomorphisms and $\epsilon$-homomorphisms, we use both to describe subnetworks and similar networks. An homomorphism transform a network to another in such a way the discrete structure giving by the first network can  lives in part of the other one, or these two networks are very similar but no equals, in particular in the probabilistic way. An $\epsilon$-homomorphism is the same but  with the condition that the probability distributions of the networks are close, and we use a preestablishes $0<\epsilon<1 $  as a distance between the probabilities.
\section{Finite dynamical systems and probabilistic Boolean networks}
Two finite dynamical systems $(X,f)$ and $(Y,g)$ are isomorphic
(or equivalents) if there exists a bijection $\phi: X\rightarrow
Y$ such that $\phi\circ f=g\circ \phi$, ( or $ f=\phi^{-1}\circ
g\circ \phi)$. If $\phi$ is not a bijection map then $\phi$ is an
homomorphism.

If $Y\subset X$ is such that $f(Y)\subset Y$ then $(Y,f|_y)$ is a
sub-FDS of $(X,f)$, where $f|_Y$ is the map restricted to Y. There exists naturally
  an injective morphism from $Y$ to $X$ called inclusion and denoted by
  $\iota$. The state space of a FDS $(X,f)$, is a digraph whit vertices the
  set $X$, and with an arrow from $u$ to $v$ if $f(u)=v$.

\medskip
\par\noindent
\textbf{Example of FDS-homomorphism} The  FDSs
$X=(\{0,1\}^2,f_1(x,y)=(xy,y))$, and
$Y=(\{0,1\}^2,f_2(x,y)=(x,(x+1)y))$ are isomorphics, because their
state spaces are isomorphics.
\[\begin{array}{cccccc}
 (1,0)&\rightarrow ^{f_1}  & (0,0)&  & (0,1)&\\
    & &\circlearrowright & &\circlearrowleft ^{f_1} & \curvearrowright^{f_1} \\
 & & & & &(1,1) \cr
\end{array} \]
\[\begin{array}{cccccc}
 (1,1)&\rightarrow ^{f_2}  & (1,0)&  & (0,0)&\\
    & &  \circlearrowright & &\circlearrowleft ^{f_2} & \curvearrowright^{f_2} \\
 & & & & &(0,1) \cr
\end{array} \]
In fact, the isomorphism $\phi:\{0,1\}^2\rightarrow \{0,1\}^2$ is  the
bijection $\phi(1,0)=(1,1)$, $\phi(0,0)=(1,0)$, $\phi(0,1)=(0,0)$,
and $\phi(1,1)=(0,1)$. The following is an example of homomorphism
(inclusion) with $Z=\{\{(0,0),(1,0)\},f_1\}$.
\[\begin{array}{ccccccccc}
(1,0)&&  & (1,0)& & &  & (0,1)&\\
  \downarrow f_1&&\hookrightarrow ^\iota & \downarrow f_1 & &  & &\circlearrowleft ^{f_1} & \curvearrowright^{f_1} \\
 (0,0)&\circlearrowleft&&(0,0)& \circlearrowleft& & & &(1,1) \cr
\end{array} \]
A Probabilistic Boolean Network $\mathcal{A}=(V,F,C)$ is defined by
the following sort (type) of objects \cite{SDZ,SDKZ}:
 a set of nodes (genes) $V=\{x_1,\ldots ,x_n\}$, $x_i\in
\{0,1\}$, for all $i$; a family  $F=\{F_1, F_2,\ldots,F_n\}$ of
ordered sets $
F_i=\{f_1^{(i)},f_2^{(i)},\ldots,f_{\ell(i)}^{(i)}\}$ of Boolean
functions $f_j^{(i)}:\{0,1\}^n\rightarrow\{0,1\}$, for all  $j$
called predictors;  and a list $C=(C_1, \ldots ,C_n)$,
$C_i=\{c_1^{(i)},\ldots ,c_{\ell(i)}^{(i)}\}$, of selection
probabilities. The selection probability that the function
$f_j^{(i)}$ is used for the vertex $i$ is
$c_j^{(i)}=Pr\{f^{(i)}=f_j^{(i)}\}$. The dynamic of the PBN is
given  by a vector of functions $\mathbf{f}
_k=(f_{k_1}^{(1)},f_{k_2}^{(2)},\ldots,f_{k_n}^{(n)})$ for
$1\le{k_i}\le{ l(i)},$ and $f_{k_i}^{(i)}\in{F_i}$, where
$k=[k_1,\ldots , k_n]$, $1\leq k_i\leq \l (i).$ The map
 $\mathbf{f}_k: {\{0,1\}}^n\rightarrow{\{0,1\}} ^n$
acts as a transition function. Each variable $x_i\in\{0,1\}^n$
represents the state of the vertex $i$.  All functions are updated
synchronously. At every time step, one of the functions is
selected randomly from the set $F_i$ according to a predefined
probability distribution. The selection probability that the
transition function  $\mathbf{f}
_k=(f_{k_1}^{(1)},f_{k_2}^{(2)},\ldots,f_{k_n}^{(n)})$ is used to
 go  from the state $u\in \{0,1\}$ to another state $\mathbf{f}
_k(u)=v\in \{0,1\}^n$ is given by
 \[c_{\mathbf{f}_k}=\prod_{i=1}^n c_{k_i}^{(i)}.
\]
The dynamical transition structure of a PBN can be described by a
Markov chain with fixed transition probabilities. There are two
digraphs structures associated with a PBN: the low-level digraph
$\Gamma$, consisting of genes functions essentiality relations;
and the high-level digraph which consists of the states of the
system and the transitions between states. The  matrix $T$
associated to the high level digraph formed by placing $p(u,v)$ in
row $u$ and column $v$, where $u,\ v\  \in \{0,1\}^n$ is called
the transition probability matrix or chain matrix,
$p(u,v)=\sum_{\mathbf{f_k| f_k}(u)=v}c_{\mathbf{f}_k}$.

\section{Probabilistic Regulatory Gene Networks}
 A Probabilistic Gene Regulatory Network (PRN) is
a triple $\mathcal{X}=(X,F,C)$ where $X$ is a finite set and
$F=\{f_1, \ldots , f_n\}$ is a set of functions from $X$ into
itself, with a list $C=(c_1, \ldots ,c_n)$ of selection
probabilities, where $c_i=p(f_i)$.
 We associate with each PRN  a weighted digraph, whose vertices are the elements of $X$, and
if $u,v\in X$, there is an arrow going from $u$ to $v$ for each
function $f_i$ such that $f_i(u)=v$, and the probability $c_i$ is
assigned to this arrow. This weighted digraph will be called the state
space of $\mathcal{X}$. In this paper, we use the notation PRN for one or more networks.\\
\begin{example}
\end{example}
If $X=\{0,1\}^2$,
$F=\{f_1(x,y)=(x,y), f_2(x,y)=(x,0),$\\$f_3(x,y)=(1,y),f_4(x,y)=(1,0)\}$; and $C=\{.46,.21,.22,.11\}$,   the state space  of
$\overline{\mathcal{X}}=(X,F,C)$ is  the following:\\
\[\begin{array}{c}
\overset{.67}{\circlearrowright}(0,0)\leftarrow ^{.21} (0,1)\circlearrowleft^{.46}\\
^{.33}\downarrow \hbox{   }\hbox{  }\swarrow_{ .11}\downarrow^{.22}\\
\overset{1}{\circlearrowright}(1,0)\overset{.32}{\longleftarrow} (1,1)\circlearrowleft^{.68} \cr
\end{array} T=\left[\begin{array}{cccc} .67&0&.33&0\\
.21&.46&.11&.22\\
0&0&1&0\\
0&0&.32&.68\cr
\end{array}\right]\]
\section{Homomorphisms and $\epsilon$-homomorphisms of PRN}
If $C$ is a set of selection probabilities we denote by $\chi$ the
characteristic function over $C$. That is $\chi:C\cup\{0\}
\rightarrow \{0,1\}$ such that $\chi(c)=1$, if $c\ne 0$ and
$\chi(0)=0$. Let
$\mathcal{X}_1=(X_1,F=(f_i)_{i=1}^n ,C)$ and
$\mathcal{X}_2=(X_2,G=(g_j)_{j=1}^m, D)$ be two PRN.\\
\begin{definition}[ Homomorphisms of PRN] A  map $\phi:X_1\rightarrow X_2$ is  an
\textbf{homomorphism} from $\mathcal{X}_1$ to $
\mathcal{X}_2$, if for all $f_i$ there exists a $g_j$, such that
for all $u$, $v$ in $\mathcal{X}_1$,
\[\hbox{(1) } \phi \circ f_i=g_j\circ \phi; \hbox{ and  } \hbox{(2) }\chi(d_{g_j}(\phi(u),\phi(v)))\geq
\chi(c_{f_i}(u,v)).\]
\[\begin{array}{c}
X_1\overset{f_i}{\longrightarrow} X_1\\
\phi \downarrow  \hspace{.4in} \downarrow \phi \\
X_2\overset{g_j}{\longrightarrow} X_2\cr
\end{array}\]
If $\phi:X_1\rightarrow X_2$ is a bijective map, then $\phi$  is an isomorphism .\end{definition}
  \begin{example} [PRN-Homomorphism]
  \end{example}
If    $\overline{\mathcal{X}}=(X;F;C)$ is the PRN in Example $1$, and  $\mathcal{X}_1=(X;F'=\{f_1, f_2,f_3\};C'=\{.47,.28,.25\})$ is a new PRN over the same set $X$ with  different probabilities and only three functions.
\[\begin{array}{c}
\mathcal{X}_1\\
^{.75}\circlearrowright(0,0)\leftarrow ^{.28} (0,1)\overset{.47}{\circlearrowleft}\\
^{.25}\downarrow \hspace{.4in}\hbox{  }\downarrow^{.25}\\
\circlearrowright^{1}(1,0)\overset{.28}{\longleftarrow }(1,1)\circlearrowleft^{.72} \cr
\end{array}  \overset{\phi}{\hookrightarrow}    \begin{array}{c}
\overline{\mathcal{X}}\\
^{.67}\circlearrowright(0,0)\leftarrow ^{.21} (0,1)\overset{.46}{\circlearrowleft}\\
^{.33}\downarrow \hbox{   }\hbox{  }\swarrow_{ .11}\downarrow^{.22}\\
\circlearrowright^{1}(1,0)\overset{.32}{\longleftarrow} (1,1)\circlearrowleft^{.68} \cr
\end{array}\]
\[T_1=\left[\begin{array}{cccc} .75&0&.25&0\\
.28&.47&0&.25\\
0&0&1&0\\
0&0&.28&.72\cr
\end{array}\right],\overline{T}=\left[\begin{array}{cccc} .67&0&.33&0\\
.21&.46&.11&.22\\
0&0&1&0\\
0&0&.32&.68\cr
\end{array}\right]\]
The homomorphism $\phi : \mathcal{X}_1\rightarrow \overline{\mathcal{X}}$ is a bijective map, $\phi (x)=x$, over the set of states, but an inclusion over the set of arrows, because the arrow going from  $(0,1)$ to $(1,1)$ in $\overline{\mathcal{X}}$ doesn't appear in  $\mathcal{X}_1$.
 The first condition for homomorphism is obvious. The condition (2) holds, because the inclusion of arrows. The two transition matrices are connected by this inclusion, since  if the place $ij$ in the first matrix $\ne 0$ then this place is $\ne 0$  in the second network too. The two PRN are not isomorphics because the probabilities are not equals. Since, there are no specific condition about the probability distribution in both PRN, we include a third condition, obtaining in this way a new  concept that we will call $\epsilon$-homomorphism of PRN.

\textbf{Condition (3) for $\epsilon$-Homomorphism}\ \emph{The
distributions of probabilities following the homomorphism are
enough close}. \emph{An $\epsilon$- homomorphism is an
homomorphism  that satisfies the condition, for all $i$, $j$, $max|p(u_i,u_j)-p(\phi(u_i),\phi(u_j))|\le \epsilon$, where $\epsilon >0$ is  a real number that we previously  determine for the  applications.}

As a consequence of this condition, if we use a test, as
 Kolmogorov-Smirnov test, the differences between the two distributions are  $\le \epsilon$ again. In order to determine $\epsilon$ for the homomorphism, we use the transition matrices. In the above example $\epsilon =.11$.\\
\[T_1-\overline{T}=\left[\begin{array}{cccc} .08&0&-.08&0\\
.07&.01&-.11&.03\\
0&0&0&0\\
0&0&-.04&.04\cr
\end{array}\right]\]
\textbf{Conclusion}

If the homomorphism is a bijective map like here, the transition matrices $T_1$ and $T_2$ have the same order, and   $\sum_{i=1}^n (T_1-T_2)_{ij}=0$, for $j=\bar{1,n}$\\
\begin{theorem}
If $\phi:\mathcal{X}_1\rightarrow\mathcal{X}_2$ is an $\epsilon$-homomorphism, then the transition matrices $T_1$ and $T_\phi$ satisfy the condition: \\$max|({T_1}^n)_{ij}-({T_\phi}^n)_{ij}|\leq \epsilon$,
 for all possible $i$ and $j$, and all $n>1$. If the homomorphism is injective and $\epsilon< 1 $, the steady state of $T_1$ and the steady state of $T_\phi$ are close, that is satisfy \\
  $\begin{array}{c}|\pi_1-\pi_\phi|=max_i|\pi_1(i)-\pi_\phi(i)|\leq\epsilon\cr
 \end{array}$.\end{theorem}
 \begin{proof}
 It is clear if  we do the following
\[|p(u,f^2(u))-p(\phi(u),\phi(f^2(u))|=\]
\[|p(u,f(u))p(f(u),f^2(u))-p(\phi(u),\phi(f(u)))p(\phi(f(u)),\phi(f^2(u)))|\]
\[\leq |p(f(u),f^2(u))||p(u,f(u))-p(\phi(u),\phi(f(u)))|+\]
\[|p(\phi(u),\phi(f(u)))||p(f(u),f^2(u))-p(\phi(f(u)),\phi(f^2(u)))|\leq\epsilon\]
Then our  aim holds.
 \end{proof}
 \section{Algebra of Probabilistic Regulatory Networks}
 \medskip
\par\noindent
\textbf{Sum of two PRN}

 Let $\mathcal{X}_1=(X_1,F=(f_i)_{i=1}^n ,C)$ and
$\mathcal{X}_2=(X_2,G=(g_j)_{j=1}^m, D)$ be two PRN. The sum
$\mathcal{X}_1\oplus\mathcal{X}_2=(X_1\dot{\cup} X_2, F\vee G,C
\vee D)$ is a  PRN where
\begin{itemize}
\item  [(1)]$X_1\dot{\cup} X_2$ is the disjoint union of $X_1$ and
$X_2$.
 \item [(2)]the
function $h_{ij}=(f_i\vee g_j)$ is defined by $h_{ij}(x)=f_i(x)$
if $x\in X_1$ and $h_{ij}(x)=g_j(x)$ if $x\in X_2$.

\item [(3)]the probability $p(h_{ij})=c_i \vee d_j$, that is
$p(h_{ij})=c_i$ if $h_{ij}=f_i$ or $p(h_{ij})=d_j$ if
$h_{ij}=g_j$.
\end{itemize}

If $T_1$ and $T_2$ are the transition matrices of $\mathcal{X}_1$ and $\mathcal{X}_2$ respectively,
Then $T=\left( \begin{array}{cc}
T_1 &0\\
0& T_2\cr
\end{array}\right)$ is the transition matrix of $\mathcal{X}_1\oplus\mathcal{X}_2$.
\begin{example}
\end{example} An example of sum is the PRN obtained by summing
the same PRN twice, $\mathcal{X}\oplus \mathcal{X}$. To make the
disjoint union, we  subindicate  $X$ with $0$ for the first $X$
and with $1$ for the second $X$. That is, the new set is
\[X_0\dot{\cup} X_1=\{(0,0,0),(0,1,0),(0,1,0),(1,1,0)\}\]
\[\cup \{(0,0,1),(0,1,1),(0,1,1),(1,1,1)\}.\]
The digraph is:
\[\begin{array}{ccccc}
   \curvearrowright^{.6} & & \curvearrowright ^{.4} & &\curvearrowright ^{1} \\
 (1,1,0)&\rightarrow ^{.4} & (1,0,0)& \rightarrow ^{.6} & (0,0,0)\\
    & &  & &  \curvearrowright ^{1}\\
 & & & & (0,1,0) \cr
\end{array}\]
\[\begin{array}{ccccc}
   \curvearrowright^{.6} & & \curvearrowright ^{.4} & &\curvearrowright ^{1} \\
 (1,1,1)&\rightarrow ^{.4} & (1,0,1)& \rightarrow ^{.6} & (0,0,1)\\
    & &  & &  \curvearrowright ^{1}\\
 & & & & (0,1,1) \cr
\end{array}  \]

This is a way to construct a PRN over $\{0,1\}^n$ using either one
or two  PRN over $\{0,1\}^{n-1}$, since $2^{n-1}+2^{n-1}=2^n$.
\medskip
\par\noindent
\textbf{Superposition}

It is clear that a PRN is the
\textbf{superposition} of several Finite dynamical Systems (FDS)\cite{H} over the same set $X$ with
probabilities assigned to each FDS. Since each functions defined
over a finite field can be wrote as a polynomial function, we will
use this notation for functions over a finite field, \cite{AGM}. If
$X=\{0,1\}=\Z_2$, the finite field of two elements,  all the FDSs
over $X$ have one of the following state space, where $f_1(x)=x;\
f_2(x)=1;\ f_3(x)=0;\ f_4(x)=x+1$, $\forall x\in X$:

\[\begin{array}{c}
L_1\\
   0\circlearrowleft\\
 \\
1\circlearrowleft\cr
\end{array}
\begin{array}{c}
L_2\\
0\\
 \downarrow\\
1\circlearrowleft\cr
\end{array}
\begin{array}{c}
L_3\\
0\circlearrowleft\\
 \uparrow\\
 1\cr
\end{array}  \begin{array}{c}
L_4\\
0\\
\uparrow \downarrow\\
1\cr

\end{array}  \]
If $p_i$ denotes de probability assigned to $L_i$, and $T_i$ denotes its transition matrix,
 then the set of all PRN is described as follows.
\[\left\{ (X,F,C)|T=\sum _{i=1}^4 p_iT_i=\left(\begin{array}{cc}
   p_1+p_3&p_2+p_4\\
p_3+p_4& p_1+p_2\cr
\end{array} \right) \sum _{i=1}^4 p_i=1
\right\}\]
We denote by $L_1L_2$ the superposition of $L_1$ and $L_2$, and similarly  $L_1L_3$ is the superposition of $L_1$ and $L_3$. The state spaces are the following:
 \[\begin{array}{c}
 L_1L_2\\
   0\circlearrowleft^{p_1}\\
 \downarrow ^{p_2}\\
1\circlearrowleft^{1}\cr
\end{array},  \  \   \begin{array}{c}
L_1L_3\\
0\circlearrowleft^{1}\\
 \uparrow^{p_3}\\
 1 \circlearrowleft^{p_1}\cr
\end{array} , \ \begin{array}{c}
L_1L_4\\
0\circlearrowleft^{p_1}\\
 \overset{p_4}{\uparrow \downarrow}\\
 1 \circlearrowleft^{p_1}\cr
\end{array} ,  \]
\[\begin{array}{c}
L_2L_3\\
   0\circlearrowleft^{p_3}\\
 \downarrow ^{p_2}\uparrow^{p_3}\\
1\circlearrowleft^{p_3}\cr
\end{array}  \  \   \begin{array}{c}
L_2L_4\\
\hspace{-.3in}0\\
 \downarrow^{p_4} \uparrow  ^{p_2+p_4}\\
 1 \circlearrowleft^{p_3}\cr
\end{array} , \]
\[\begin{array}{c}
L_3L_4\\
0\circlearrowleft^{p_3}\\
 \overset{p_4}{\downarrow}\uparrow^ {p_4+p_3}\\
\hspace{-.2in} 1 \cr
\end{array},  \]

For example, with transition matrices \[T_{12}=T_1+T_2=\left(\begin{array}{cc}
   p_1&p_2\\
0& 1\cr
\end{array} \right) \  T_{13}=T_1+T_3=\left(\begin{array}{cc}
   1&0\\
p_3& p_1\cr
\end{array} \right)\]
\medskip
\par\noindent
\textbf{Product of two PRN}

Let
$\mathcal{X}_1=(X_1,F=(f_i)_{i=1}^n ,C)$ and
$\mathcal{X}_2=(X_2,G=(g_j)_{j=1}^m, D)$ be two PRN. The product
$\mathcal{X}_1\times \mathcal{X}_2=(X_1\times X_2, F \times G,C
\wedge D)$ is a  PRN where
\begin{itemize}
\item  [(1)]$X_1\times X_2$ is the cartesian product  of $X_1$ and
$X_2$.
 \item [(2)]the
function $h_{ij}=(f_i,g_j)$ is defined by
\[h_{ij}(x_1,x_2)=(f_i(x_1),g_j(x_2))\]  for $x_1\in X_1$, and
$x_2\in X_2$.

\item [(3)]the probability $p(h_{ij})$ is a function of $c_i$ and $d_j$, for example $p(h_{ij})=\frac{c_i+d_j}{2}$.
\end{itemize}
\begin{example}
\end{example}
The product $L_1L_2\times L_1L_3$ is the PRN with four states $\{(0,0),(0,1),$
$(1,0),(1,1)\}$ and
four functions \[f_{11}(x,y)=(x,y), f_{13}(x,y)=(x,0),\]
 \[ f_{21}(x,y)=(1,y),\  f_{23}(x,y)=(1,0).\]
The state space is the following:
\[\begin{array}{c}
L_1L_2\times L_1L_3\\
^{p_{11}+p_{13}}\circlearrowright(0,0)\leftarrow ^{p_{13}} (0,1)\circlearrowleft^{p_{11}}\\
^{p_{23}+p_{21}}\downarrow \hbox{   }\hbox{  }\swarrow_{ p_{23}}\downarrow^{p_{21}}\\
\hbox{      }\hbox{      }\hbox{ }\hbox{  }\hbox{   }\hspace{.3in}    \hbox{     } \  \hbox{    }  \circlearrowright^{1}(1,0)\longleftarrow (1,1)\circlearrowleft^{p_{11}+p_{21}} \\
\hspace{.3in} ^{p_{13}+p_{23}} \cr
\end{array}\]
The transition matrix is the following
\[T=\left(\begin{array}{cccc}
 p_{11}+p_{13}  &0&p_{23}+p_{21}&0\\
p_{13}& p_{11}&p_{23}&p_{21}\\
0&0&1&0\\
0&0&p_{13}+p_{23}& p_{11}+p_{21}\cr
\end{array} \right)\]
\subsection{Linear Probabilistic Regulatory Networks }

 A linear PRN is a
\textbf{\emph{superposition}} of linear FDS. A linear FDS is a
pair $(X,f)$ where $f$ is a linear function, and $X$ is a vector
space over a finite field. So, a linear PRN is a triple $(X,(f_i)_{i=1}^m, C)$, where $X$ is
a finite vector space, the functions $f_i:X\rightarrow X$ are linear functions, and
 $C=\{c_i=p(f_i)\}$. The set $X$ has cardinality a power of a prime number and each linear function
 is determined by its characteristic polynomial and the companion matrix.

  If $X=\Z_3=\{ 0,1,2\}$ is the field  of
integer modulo $3$, then the linear functions are:  $f_1(x)=x$,
$f_2(x)=2x$, and $f_3(x)=0$ for all $x\in \Z_3$. So, the  linear
PRN are the following:
\[\hbox{   }\begin{array}{c}
\{ f_1,f_2\}\\
  \circlearrowleft ^{1}\\
   \mathbf{0} \\
   \\
  \circlearrowright^{p_1} \mathbf{1}\rightleftarrows ^{p_2}\mathbf{2}\circlearrowleft^{p_1}\cr
\end{array} \hbox{   }
\hbox{    }\begin{array}{c}
\{f_1,f_3\}\\
\circlearrowleft ^{1}\\
   \mathbf{0} \\
\nearrow^{p_3}\nwarrow  \\
\mathbf{1} \circlearrowleft _{p_1}
 \hbox{    }\circlearrowright_{p_1} \mathbf{2}
\cr
\end{array}\]
\[\begin{array}{c}
\{f_2,f_3\}\\
\circlearrowleft ^{1}\\
   \mathbf{0} \\
\nearrow^{p_3}\nwarrow  \\
\mathbf{1} \rightleftarrows _{p_2}  \mathbf{2}
\cr
\end{array}
\hbox{    }\begin{array}{c}
\{f_1,f_2,f_3\}\\
\circlearrowleft ^{1}\\
   \mathbf{0} \\
\nearrow^{p_3}\nwarrow  \\
\circlearrowleft _{p_1}\mathbf{1} \rightleftarrows _{p_2}  \mathbf{2}\circlearrowleft _{p_1}
\cr
\end{array}\]
If $X=\Z_2\times \Z_2$ is the vector space with $4$ elements over
the field $\Z_2$, then there are $4$ linear FDS not isomorphics. In fact,
using matrix, the possible characteristics polynomials $p_f(\lambda)$ are:
$\lambda ^2,\ \lambda ^2+\lambda \ \lambda ^2 +1,\ \lambda
^2+\lambda +1$. The companion matrices of these linear functions are:
\[A_1\left(\begin{array}{cc}
0&0\\
0 &0\cr
\end{array}\right) \hbox{   }A_2\left(\begin{array}{cc}
0&0\\
0 &1\cr
\end{array}\right) \  A_3\left(\begin{array}{cc}
1&0\\
0 &1\cr
\end{array}\right)\hbox{   }A_4\left(\begin{array}{cc}
0&1\\
1 &1\cr
\end{array}\right)\]
Then the FDS associated to this matrices are:
\[\begin{array}{ccc}
A_1&&\\
\circlearrowright(0,0)&\leftarrow& (1,0)\\
\uparrow&\nwarrow&\\
(0,1)&&(1,1)\cr
\end{array}
\hbox{   }\begin{array}{ccc}
A_2&&\\
\circlearrowright(0,0)&\leftarrow& (1,0)\\
&&\\
\circlearrowright(0,1)&\leftarrow&(1,1)\cr
\end{array}\hbox{   }\]
\[\begin{array}{ccc}
A_3&&\\
\circlearrowright(0,0)&& \circlearrowright(1,0)\\
&&\\
\circlearrowright(0,1)&&\circlearrowright(1,1)\cr
\end{array}\hbox{    }\begin{array}{ccc}
A_4&&\\
\circlearrowright(0,0)&& (1,0)\\
&\swarrow&\uparrow\\
(0,1)&\rightarrow&(1,1)\cr
\end{array}\hbox{   }\]
The linear PRN  with two functions are the following:
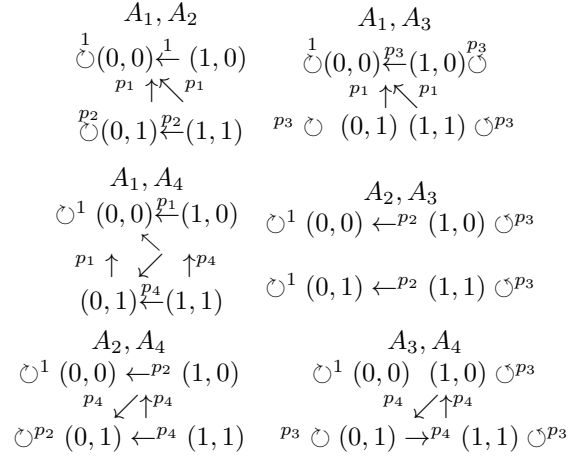
\[\begin{array}{c}
A_1,A_2\\
{\overset{1}{\circlearrowright}}(0,0){\overset{1}{\leftarrow}(1,0)} \\
^{p_1}\uparrow \nwarrow ^{p_1}\\
{\overset{p_2}{\circlearrowright}}(0,1){\overset{p_2}{\leftarrow}}(1,1)\cr
\end{array}
\begin{array}{c}
A_1,A_3\\
{\overset{1}{\circlearrowright}}(0,0){\overset{p_3}{\leftarrow}}(1,0){\overset{p_3}{\circlearrowleft}} \\
^{p_1}\uparrow\nwarrow^ {p_1}\\
^{p_3} \circlearrowright\hbox{  }(0,1) \   (1,1)\circlearrowleft ^{p_3}\cr
\end{array} \]
\[ \begin{array}{c}
A_1,A_4\\
\circlearrowright^1(0,0){\overset{p_1}{\leftarrow}}(1,0)\\
^{p_1}\uparrow \hspace{.1in}  \overset{\nwarrow}{\swarrow}\hspace{.1in} \uparrow^{p_4}\\
\  (0,1)  {\overset{p_4}{   \leftarrow}}(1,1)
\cr
\end{array}
\begin{array}{c}
A_2,A_3\\
\circlearrowright^{1}(0,0)\leftarrow^{p_2} (1,0)\circlearrowleft ^{p_3}\\
\\
\circlearrowright^{1}(0,1)\leftarrow^{p_2}(1,1)\circlearrowleft ^{p_3}\cr
\end{array}\]
\[\begin{array}{c}
A_2,A_4\\
\circlearrowright^1(0,0)\leftarrow ^{p_2} (1,0)\\
^{p_4}\swarrow \uparrow ^{p_4}\\
\circlearrowright^{p_2}(0,1)\leftarrow ^{p_4}(1,1)\cr
\end{array}\hbox{    }
\begin{array}{c}
A_3,A_4\\
\circlearrowright^1 (0,0)\  \   (1,0)\circlearrowleft ^{p_3}\\
\    ^{p_4}\swarrow\uparrow ^{p_4} \\
^{p_3}\circlearrowright(0,1)\rightarrow^{p_4}(1,1)\circlearrowleft^{p_3}
\cr
\end{array}\hbox{   }
\]

\section{Invariant Subnetworks and Projections}
A subnetwork  $Y\subseteq X$  of $\mathcal{X}=(X,F,C)$ is an
\textbf{invariant subnetwork or a sub-PRN} of $\mathcal{X}$ if
 $f_i(u)\in Y$ for all $u\in Y$, and  $f_i\in F$. Sub-PRNs
 are sections of a PRN, where there aren't arrows going out.
  The complete network $X$, and any cyclic state with
  probability 1, are sub-PRNs. An invariant subnetwork is irreducible if doesn't have a
proper invariant subnetwork.
 \emph{An endomorphism is a projection if $\pi^2=\pi$. }\\
\begin{theorem}
If there  exists a projection from $\mathcal{X}$ to a subnetwork $\mathcal{Y}$ then $\mathcal{Y}$ is  an invariant subnetwork of $\mathcal{X}$.
\end{theorem}
\begin{proof}
Suppose that there exists a projection $\pi:X\rightarrow Y$. If  $y\in Y$,  by definition of projection $\pi(y)=y$,  and $f_i(\pi(y))=\pi(g_j(y))$. Therefore all  arrows in the subnetwork $Y$ are going inside $Y$, and the network is invariant.
\end{proof}

\begin{example}
\end{example}
The PRN $\overline{\mathcal{X}}$ has two invariant subnetworks with projections \\
 $\pi_1(x,y)=(x,0)$ and $\pi_2(x,y)=(1,y)$.
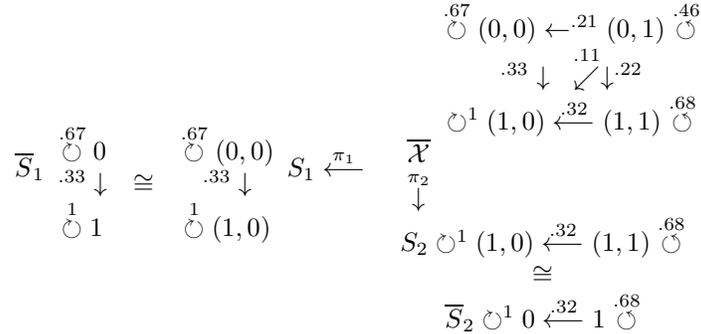
\[\overline{S}_1\begin{array}{ccc}
&&\hspace{.4in}\\
\overset{.67}{\circlearrowright}0 &&\overset{.67}{\circlearrowright}(0,0)\\
^{.33}\downarrow &\cong &^{.33}\downarrow   \\
\overset{1}{\circlearrowright} 1& &\overset{1}{\circlearrowright}(1,0)\cr
\end{array} S_1\overset{\pi_1}{\longleftarrow}\begin{array}{c}
 \hspace{.3in}\overset{.67}{\circlearrowright}(0,0)\leftarrow ^{.21} (0,1)\overset{.46}{\circlearrowleft}\\
\hspace{.3in}^{.33}\downarrow \hbox{   }\overset{.11}{\swarrow}\downarrow^{.22}\\
\hspace{.3in}\circlearrowright^{1}(1,0)\overset{.32}{\longleftarrow }(1,1)\overset{.68}{\circlearrowleft} \\
\hspace{-1.3in}\overline{\mathcal{X}}\\
\hspace{-1.3in}\overset{\pi_2}{\downarrow}\\
S_2\circlearrowright^{1}(1,0)\overset{.32}{\longleftarrow} (1,1)\overset{.68}{\circlearrowleft}\hspace{-.3in}\\
\cong\hspace{-.3in} \\
\overline{S}_2\circlearrowright^{1}0\overset{.32}{\longleftarrow} 1\overset{.68}{\circlearrowleft}\hspace{-.3in}\cr
\end{array} \]
Checking the probabilities for $\pi_1$ and $\pi_2$, we have
$\epsilon_1=.68$; and  $\epsilon_2=.67$.

We can observe that
$\overline{\mathcal{X}}\cong \overline{S}_1\times \overline{S}_2$.
\begin{example}
\end{example} The subnetwork $\mathcal{X}_1=(\{(x,y,1)\}, F, C)$ is an invariant subnetwork of $\mathcal{X}=(\{0,1\}^3,F, C) $.
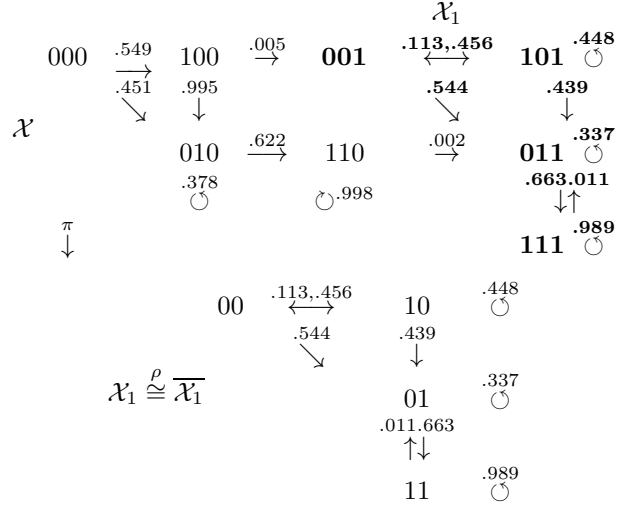
\[\mathcal{X}\begin{array}{ccccccc}
&&&&&\mathcal{X}_1&\\
000&\underrightarrow{^{.549}}&100&\overset{.005}{\rightarrow }&\textbf{001}&\overset{\textbf{.113,.456}}{\longleftrightarrow}&\textbf{101}\overset{\textbf{.448}}{\circlearrowleft} \\
\hspace{.2in}&\overset{.451}{\searrow }&\overset{.995}{\downarrow }&&&\overset{\textbf{.544}}{\searrow }&\overset{\textbf{.439}}{\downarrow} \\
&&010&\overset{.622}{\longrightarrow }&110&\overset{.002}{\rightarrow }&\textbf{011}\overset{\textbf{.337}}{\circlearrowleft }\\
&&\overset{.378}{\circlearrowleft}&&\circlearrowright ^{.998}&&\overset{\textbf{.663.011}}{\downarrow \uparrow }\\
\overset{\pi}{\downarrow}&&&&&& \textbf{111}\overset{\textbf{.989}}{\circlearrowleft }\cr
\end{array}\]
\[\mathcal{X}_1 \overset{\rho}{\cong} \overline{\mathcal{X}_1}\begin{array}{cccc}
00&\overset{.113, .456}{\longleftrightarrow}&10&\overset{.448}{\circlearrowleft }\\
&\overset{.544}{\searrow }&\overset{.439}{\downarrow }&\\
&&{01}&\overset{.337}{\circlearrowleft }\\
&&\overset{.011 .663}{\uparrow \downarrow }&\\
&& 11&\overset{.989}{\circlearrowleft }\cr
\end{array}\]
Ordering the elements in the following way $\{(0,0,0),(0,1,0),\\(1,0,0),(1,1,0),(0,0,1),(0,1,1),(1,0,1),(1,1,1)\}$, the matrix
 $ T_{X_1}=\left[\begin{array}{cccc} 0&.544&.456&0\\
0&.337&0&.663\\
.113&.448&.439&0\\
0&.011&0&.989\cr
\end{array}
\right]$ is an invariant part of the transition matrix
 $T_{\mathcal{X}}=\left[\begin{array}{cc}
T_{11}& T_{12}\\
0&T_{X_1}\cr
\end{array}\right].$\\
Using the projection $\pi:\mathcal{X}\rightarrow \mathcal{X}_1$,
$\pi(x,y,z)=(x,y,1)$; and the isomorphism $\rho (x,y,1)=(x,y)$,
the network $\mathcal{X}$ is projected over the network
$\overline{\mathcal{X}_1}$. Checking the arrows the projection
$\overline{\pi}$ is a  $.5$-homomorphism.
\subsection{Mathematical background}
\begin{theorem}
  If $\phi_1:\mathcal{X}_1\rightarrow \mathcal{X}_2$ is an
$\epsilon_1$-homomorphism, and
$\phi_2:\mathcal{X}_2\rightarrow \mathcal{X}_3$ is another $\epsilon_2$-homomorphism. Then $\phi=\phi_2\circ \phi_1:\mathcal{X}_1\rightarrow \mathcal{X}_3$ is an $\epsilon$-homomorphism. Therefore the Probabilistic Regulatory Networks with the homomorphisms of PRN form the category \textbf{PRN}.
\end{theorem}
\begin{proof}
The Probabilistic Regulatory Networks with the PRN homomorphisms  is a
category if: the composition is an homomorphism, and satisfy the
associativity law; and there exists an identity homomorphism for
each PRN.

  (1) Let $\phi_1:\mathcal{X}_1\rightarrow \mathcal{X}_2$ be an $\epsilon_1$-homomorphism, and let
$\phi_2:\mathcal{X}_2\rightarrow \mathcal{X}_3$ be an $\epsilon_2$-homomorphism.
If $q_t$, $g_k$ and $f_j$ are
functions in each PRN, and such that $\phi_1 \circ f_j=g_k \circ
\phi_1$ and $\phi_2 \circ g_k=q_t \circ \phi_2$, then we
will  prove that: $\phi \circ f_j=q_t \circ \phi.$ In fact,
 \[(\phi_2 \circ \phi_1) \circ f_j=\phi_2 \circ (\phi_1 \circ f_j)=\]
 \[\phi_2 \circ (g_k
 \circ \phi_1)=(\phi_2 \circ  g_k)\circ \phi_1=\]
 \[(q_t \circ \phi_2)\circ \phi_1=q_t \circ (\phi_2\circ \phi_1).\]
\noindent
 (2)We want to prove that
\[
\chi(t_k(\phi(u),\phi(v)))\ge{\chi(c_i(u,v))}.
\]
Suppose that
${\chi(c_i(u,v))}=1$.
Then, since $\phi_1$ is an homomorphism of PRN, we have that
\[\chi(d_j (\phi_1(u),\phi_1(v)))\ge {\chi (c_i(u,v))%
}\]
 which is 1. Since
$\phi_2$ is an homomorphism of PRN, we obtain that
\[\chi(t_k(\phi(u),
\phi(v)))=\chi(t_k(\phi_2(\phi_1(u)), \phi_2(\phi_1(v))))\]

\[ \ge \chi(c_j(\phi_1(u),(\phi_1(v))=1.\]
Therefore we obtain that
\[
\chi (t_k(\phi_2(\phi_1(u)), \phi_2(\phi_1(v))))=1.
\]

Then the composition of two PRN-homomorphisms  is an homomorphism.

(3) To verify the third condition for $\epsilon$-homomorphism, we do the following.
If $p(\phi(u_1),\phi(u_2))>1$, with $u_1,\ u_2\in X_1$, then we need to prove that there exists an $\epsilon$ such that \[|p(u_1,u_2)-p(\phi(u_1),\phi(u_2))|<\epsilon.\] In fact:

\[|p(u_1,u_2)-p(\phi(u_1),\phi(u_2))|=|p(u_1,u_2)-p(\phi_1(u_1),\phi_1(u_2))+ \] \[p(\phi_1(u_1),\phi_1(u_2))-p(\phi_2(\phi_1(u_1)),\phi_2(\phi_1(u_2)))|\]
\[< |p(u_1,u_2)-p(\phi_1(u_1),\phi_1(u_2))|+\]
\[|p(\phi_1(u_1),\phi_1(u_2))-p(\phi_2(\phi_1(u_1)),\phi_2(\phi_1(u_2)))|\leq \epsilon_1+\epsilon _2\]
\[|p(u_1,u_2)-p(\phi(u_1),\phi(u_2))|< \epsilon_1+\epsilon _2\]
because $\phi_1$ and $\phi_2$ are $\epsilon$-homomorphisms.

The associativity and identity laws are easily checked, therefore
our claim holds, and \textbf{PRN} is a category.
\end{proof}
It is clear that, the PRN with the homomorphism between them form a category that we will denote $\mathcal{PRN}$. The category \textbf{PRN} is a  subcategory of $\mathcal{PRN}$,  since an homomorphism is not always an homomorphism for some $\epsilon \in \R$ enough small. But, if we don't include the condition for $\epsilon$ to be enough small, the two categories are the same, because always an homomorphism is an $\epsilon$-homomorphism for some $ \epsilon \in \R$.
\begin{theorem}\label{product}
Let $\mathcal{X}_1\times \mathcal{X}_2=(X_1\times X_2,H,E)$ be a product of PRN $\mathcal{X}_1=(X_1,F,C)$ and $\mathcal{X}_2=(X_2,G,D)$. If $\delta_i:X \rightarrow X_i$ are two PRN-homomorphisms, then there exists an homomorphism $\delta: X\rightarrow X_1\times X_2$, such that $\phi_i\circ \delta=\delta_i$ for $i=1,2$. That is, the following diagram commutes
\[\begin{array}{c}
  \hspace{.05in} X_1\times X_2 \hspace{.05in}\\
 \overset{\phi_1}{ \swarrow}\hspace{.1in}\overset{\delta}{\uparrow}\hspace{.1in}\overset{\phi_2}{\searrow}\\
  X_1\overset{\delta_1}{\longleftarrow} X \overset{\delta_2}{\longrightarrow} X_2\cr
  \end{array}\]
  This homomorphism is unique.
\end{theorem}
\begin{proof}
The function $\delta:X\rightarrow X_1\times X_2$ is defined as follows $\delta(x)=(\delta_1(x),\delta_2(x))$,  $x\in X$. $\delta$ is an homomorphism, in fact:

(1) Let $\mathcal{X}=(X,L,P)$ be a PRN. Since $\delta_1$ and $\delta_2$ are homomorphism, for all function $l_t\in L$ there exist two functions $f_i\in F$ and $g_j\in G$, such that $\delta_1\circ l_t=f_i\circ \delta_1$, and $\delta_2\circ l_t=g_j\circ \delta_2$. Then for the function $l_t$ there exists the function $(f_i,g_j)$ that satisfies $\delta \circ l_t=(f_i,g_j)\circ \delta$.
\[(\delta \circ l_t)(x)=\delta(l_t(x))=(\delta_1(l_t(x)),\delta_2(l_t(x)))=\]
\[(f_i(\delta_1(x)),g_j(\delta_2(x)))=((f_i,g_j)\circ\delta)(x)\]

(2) In order to prove $\chi(e_{ij}(\delta(x),\delta(x')))\geq \chi(p_{l_t}(x,x'))$,  suppose  $\chi(p_{l_t}(x,x'))=1$. Then $l_t(x)=x'$, and $\delta(x')=\delta(l_t(x))=(f_i,g_j)(\delta(x))$ by  part (1). Therefore $\chi(e_{ij}(\delta(x),\delta(x')))=1$, and our claim holds.

It is easy to check that $\phi_i\circ \delta=\delta_i$, in fact \[\phi_1 (\delta(x))=\phi_1(\delta_1(x),\delta_2(x))=\delta_1(x),\] for all $x\in X$.
\end{proof}
If $\delta_i$, $i=1,2$,  are $\epsilon _i$-homomorphism then \[max|p(x,x')-p(\phi_1(\delta(x)),\phi_1(\delta(x')))|\leq \epsilon_1.\] But
\[|p(x,x')-p(\delta(x),\delta(x'))+p(\delta(x),\delta(x'))-
p(\phi_1(\delta(x)),\phi_1(\delta(x')))|\leq \]
\[|p(x,x')-p(\delta(x),\delta(x'))|+\]
\[|p(\delta(x),\delta(x'))-
p(\phi_1(\delta(x)),\phi_1(\delta(x')))|\leq \epsilon_1.\]
Therefore
\[|p(x,x')-p(\delta(x),\delta(x'))|\leq \epsilon_1-|p(\delta(x),\delta(x'))-
p(\phi_1(\delta(x)),\phi_1(\delta(x')))|\]
\[|p(x,x')-p(\delta(x),\delta(x'))|\leq \epsilon_1-\overline{\epsilon}_1.\]
Therefore $\delta$ is an $\epsilon$-homomorphism. So, the theorem holds for $\epsilon$-homomorphism.

It is an immediate consequence the following result, also is true for $\epsilon$-homomorphisms.
\begin{theorem}
Let $\mathcal{X}_1\oplus \mathcal{X}_2=(X_1\times X_2,H,E)$ be a product of PRN $\mathcal{X}_1=(X_1,F,C)$ and $\mathcal{X}_2=(X_2,G,D)$. If $\gamma_i:X_i \rightarrow X$ are two PRN-homomorphisms, then there exists an homomorphism $\gamma: X_1\oplus X_2\rightarrow X $, such that $\gamma \circ \iota_i=\gamma_i$ for $i=1,2$. That is, the following diagram commutes
\[\begin{array}{c}
  \hspace{.05in} X_1\oplus X_2 \hspace{.05in}\\
 \overset{\iota_1}{ \nearrow}\hspace{.1in}\overset{\gamma}{\downarrow}\hspace{.1in}\overset{\iota_2}{\nwarrow}\\
  X_1\overset{\gamma_1}{\longrightarrow} X \overset{\gamma_2}{\longleftarrow} X_2\cr
  \end{array}\]
  This homomorphism is unique.
\end{theorem}

\begin{theorem}[Fundamental Theorem] All  reducible PRN is either a product of its non trivial sub-PRN  or a subnetwork of this product.
\end{theorem}
\begin{proof} It is trivial by definition of Product and sub-PRN.
\end{proof}
\section{Conclusions}
The intersection, and the union of two sub-PRN  is a sub-PRN, therefore the class of sub-PRN of a particular PRN is a lattice. Reduction mappings described in \cite{ID} and defined for PBN using the influence of a gene, for example $x_n$, on the predictor function $f_j^{(i)}$ to determine the selected predictor, can be extended to PRN. In order to extend this procedure to more than boolean functions, we use the polynomial description of genetic functions given in \cite{AGM}, the partial derivative is the usual in calculus  and all the concepts  in \cite{ID} can be using for PRN. Similarly our definition of projection, the reduction mappings are $\epsilon$-homomorphisms, and we can use  for genes with more than two quantization, since this extension is not a trivial work we  develop the theory and methods  in \cite{A}.

\end{document}